\documentclass[11pt,twoside,a4paper,dvips]{article}

\usepackage{amsmath,amsfonts,amssymb,amsthm}
\usepackage[all]{xy}

\swapnumbers
\theoremstyle{plain}
\newtheorem{lemma}{Lemma}[section]
\newtheorem{prop}[lemma]{Proposition}
\newtheorem{thm}[lemma]{Theorem}
\newtheorem{thmdefn}[lemma]{Theorem--Definition}
\newtheorem{cor}[lemma]{Corollary}

\newtheorem*{claim*}{Claim}

\newtheorem{nfi}[lemma]{NFI-type inequality}

\newtheorem{smptsthm}[lemma]{Smooth Points Theorem}

\theoremstyle{definition}
\newtheorem{defn}[lemma]{Definition}
\newtheorem{defns}[lemma]{Definitions}
\newtheorem{emp}[lemma]{}
\newtheorem{ex}[lemma]{Example}

\newtheorem{rk}[lemma]{Remark}

\newtheorem*{notn*}{Notation}
\newtheorem*{rk*}{Remark}
\newtheorem*{rks*}{Remarks}
\newtheorem*{term*}{Terminology}

\newtheorem{comments_on_excl_sec}[lemma]{The exclusion arguments of
\S\ref{sec:excl}}
\newtheorem{comments_on_final_sec}[lemma]{The derivation of
Theorem~\ref{thm:75main} from Theorem~\ref{thm:75aux}
in~\S\ref{sec:class}}
\newtheorem{crvs_fams34etc}[lemma]{Curves for families~34,~88 and~90}
\newtheorem{absexclsingpts}[lemma]{Absolute exclusion of singular
points with $B^3 < 0$ and $T \sim bB+cE$, $b,c>0$}
\newtheorem{condexclsingpts}[lemma]{Singular points with $B^3<0$ and
$T \sim bB$, $b>0$}

\newcommand{\CC}{\mathbb C}
 
\newcommand{\HH}{\mathcal{H}}
\newcommand{\II}{\mathcal I}
\newcommand{\LL}{\mathcal L} 

\newcommand{\NN}{\mathbb N}

\newcommand{\Oh}{\mathcal O}
\newcommand{\PP}{\mathbb{P}}
\newcommand{\PPP}{\mathcal P}
\newcommand{\QQ}{\mathbb{Q}}
\newcommand{\RR}{\mathbb R}

\newcommand{\ZZ}{\mathbb{Z}}

\newcommand{\B}{\operatorname{B}}

\newcommand{\Bs}{\operatorname{Bs}}
\newcommand{\Centre}{\operatorname{Centre}}

\newcommand{\cont}{\operatorname{cont}}
\newcommand{\CS}{\operatorname{CS}}
\newcommand{\CSXnH}{\operatorname{CS}\left(X,\textstyle\frac{1}{n}\HH\right)}

\newcommand{\dashto}{\dashrightarrow}

\newcommand{\Exc}{\operatorname{Exc}}

\newcommand{\im}{\operatorname{im}}
\newcommand{\iso}{\simeq}
\newcommand{\KXnH}{K_X + \textstyle\frac{1}{n}\HH}
\newcommand{\LCS}{\operatorname{LCS}}
\newcommand{\mult}{\operatorname{mult}}
\newcommand{\NEbar}{\operatorname{\overline{NE}}}

\newcommand{\qeq}{\sim_\QQ}

\newcommand{\VXnH}{\operatorname{V_0}\left(X,\textstyle\frac{1}{n}\HH\right)}

\newcommand{\XnepH}{\left(X,\left(\textstyle\frac{1}{n} + \ep\right)\HH\right)} 
\newcommand{\XnH}{\left(X,\textstyle\frac{1}{n}\HH\right)}

\newcommand{\ep}{\varepsilon}
\newcommand{\fie}{\varphi} 
\newcommand{\ka}{\kappa}

\begin{document}

\title{Classification of elliptic and K3 fibrations birational to some
$\QQ$-Fano 3-folds}
\author{Daniel Ryder\footnote{Email: \texttt{dj@ryder144.fsnet.co.uk}}}
\date{August 2005}
\maketitle

\begin{abstract}
{\sloppy
\noindent A complete classification is presented of elliptic and K3
fibrations birational to certain mildly singular complex Fano
3-folds.  Detailed proofs are given for one example case, namely that
of a general hypersurface $X$ of degree 30 in weighted $\PP^4$ with
weights 1,4,5,6,15; but our methods apply more generally.  For
constructing birational maps from $X$ to elliptic and K3 fibrations we
use Kawamata blowups and Mori theory to compute anticanonical rings;
to exclude other possible fibrations we make a close examination of
the strictly canonical singularities of $\XnH$, where
$\HH$ is the linear system associated to the putative birational map
and $n$ is its anticanonical degree.
}
\end{abstract}

\section{Introduction}

In \cite{CPR} Corti, Pukhlikov and Reid proved that a
general quasismooth complex variety \linebreak $X = X_d \subset
\PP(1,a_1,\ldots,a_4)$ in one of the `famous~95 families' of
$\QQ$-Fano 3-fold weighted hypersurfaces is \emph{birationally rigid}
--- that is, if $X$ is birational to some Mori fibre space $Y/S$ then
in fact $Y \iso X$.  A related problem is to classify elliptic and K3
fibrations birational to general hypersurfaces in these families; it
was Ivan Cheltsov \cite{Ch00} who first proved classification results
of this kind for several birationally rigid smooth Fano varieties,
including a general quartic 3-fold $X_4 \subset \PP^4$ and a double
cover of $\PP^3$ branched in a general sextic surface, i.e., $X = X_6
\subset \PP(1,1,1,1,3)$. In \cite{Ry02} the classification of elliptic
and K3 fibrations birational to general members of the remaining 93
families was addressed, but completed only for family~5, $X_7 \subset
\PP(1,1,1,2,3)$. In the present paper we aim firstly to give concise
proofs of some of the more generally applicable results of \cite{Ry02}
and secondly to present a complete proof of the following theorem for
family 75, which is the family referred to in the abstract.
Furthermore, we state similar theorems for families 34, 88 and 90;
these can be proved using essentially the same techniques.

\begin{thm} \label{thm:75main}
Let $X = X_{30} \subset \PP(1,4,5,6,15)_{x,y,z,t,u}$ be a general
member of family~75 of the~95.
\begin{itemize}
\item[(a)] Suppose $\Phi \colon X \dashto Z/T$ is a birational map
from $X$ to a K3 fibration $g \colon Z \to T$ (see~\ref{defns:fibr}
below for our assumptions on K3 fibrations, and also on elliptic
fibrations and Fano 3-folds). Then there exists an isomorphism
$\PP^1 \to T$ such that the diagram below commutes, where $\pi =
(x^4,y) \colon X \dashto \PP^1$.
\[ \xymatrix@C=1.6cm{
X \ar@{-->}[r]^{\Phi} \ar@{-->}[d]_{\pi} & Z \ar[d]^g \\
\PP^1 \ar[r]^{\iso} & T \\
} \]
\item[(b)] There does not exist an elliptic fibration birational to
$X$.
\item[(c)] If $\Phi \colon X \dashrightarrow Z$ is a birational map
from $X$ to a Fano 3-fold $Z$ with canonical singularities then
$\Phi$ is actually an isomorphism (so in particular $Z \iso X$ has
terminal singularities).
\end{itemize}
\end{thm}

Part (b) of this theorem was recently proved independently by Cheltsov
and Park~\cite{CP05} using somewhat different methods.  It is an
interesting result because of its relevance to the question of whether
$\QQ$-rational points of $X$ are potentially dense: a birational
elliptic fibration is one key geometric construction used to prove
potential density (see~\cite{HT00}, \cite{Ha03} and
\cite{HT01}).  The proof presented here requires close examination of
$\CSXnH$ (see~\ref{comments_on_excl_sec}), where $\HH$ is the linear
system associated to a putative birational map and $n$ is its
anticanonical degree; in~\cite{CP05} more general methods are used.
Our approach has the advantage that the other parts of
Theorem~\ref{thm:75main} follow immediately from our complete
classification of possible sets of strictly canonical centres.

The following results are analogous to Theorem~\ref{thm:75main}; they
can be proved with the same techniques, though we do not include all
the details here.  From now on we abbreviate conclusions such as that
in Theorem~\ref{thm:75main}(a) by stating that \emph{up to a
birational twist of the base,} $g \circ \Phi = (x^4,y) \colon X
\dashto \PP^1$.  The birational twist $\PP^1 \to T$ of the base is an
isomorphism in the above case because $T$ is a smooth curve.
(We assume all fibrations are morphisms of normal varieties:
see~\ref{defns:fibr}.)

\begin{thm} \label{thm:34main}
Let $X = X_{18} \subset \PP(1,1,2,6,9)_{x_0,x_1,y,z,t}$ be a general
member of family~34 of the~95.
\begin{itemize}
\item[(a)] If $\Phi \colon X \dashto Z/T$ is a birational map from $X$
to a K3 fibration $g \colon Z \to T$ then, up to a birational twist
of the base (see above for explanation), $g \circ \Phi = (x_0,x_1) \colon X
\dashto \PP^1$.
\item[(b)] Suppose $\Phi \colon X \dashto Z/T$ is a birational map from $X$
to an elliptic fibration $g \colon Z \to T$.  Then, up to a
birational twist of the base, $g \circ \Phi = (x_0,x_1,y) \colon X \dashto
\PP(1,1,2)$.
\item[(c)] If $\Phi \colon X \dashrightarrow Z$ is a birational map from $X$
to a Fano 3-fold $Z$ with canonical singularities then
$\Phi$ is actually an isomorphism (so in particular $Z \iso X$ has
terminal singularities).
\end{itemize}
\end{thm}

\begin{thm} \label{thm:88main}
Let $X = X_{42} \subset \PP(1,1,6,14,21)_{x_0,x_1,y,z,t}$ be a general
member of family~88 of the~95.  Under assumptions corresponding to
those in Theorem~\ref{thm:34main} we can conclude as follows.
\begin{itemize}
\item[(a)]  Up to a birational twist of the base, $g \circ \Phi = (x_0,x_1)
\colon X \dashto \PP^1$.
\item[(b)]  Up to a birational twist of the base, $g \circ \Phi =
(x_0,x_1,y) \colon X \dashto \PP(1,1,6)$.
\item[(c)]  $\Phi$ is actually an isomorphism.
\end{itemize}
\end{thm}

\begin{thm} \label{thm:90main}
Let $X = X_{42} \subset \PP(1,3,4,14,21)_{x,y,z,t,u}$ be a general
member of family~90 of the~95.  Under assumptions corresponding to
those in Theorem~\ref{thm:34main} we can conclude as follows.
\begin{itemize}
\item[(a)]  Up to a birational twist of the base, $g \circ \Phi = (x,y)
\colon X \dashto \PP(1,3)$.
\item[(b)]  Up to a birational twist of the base, $g \circ \Phi = (x,y,z)
\colon X \dashto \PP(1,3,4)$.
\item[(c)]  $\Phi$ is actually an isomorphism.
\end{itemize}
\end{thm}

\subsection*{Outline of paper}

Our proof of Theorem~\ref{thm:75main} has essentially three parts; this
division of the argument is closely modelled on the approach of Cheltsov to
similar problems for smooth varieties --- see e.g.~\cite{Ch00}.  In brief,
the parts are: constructing the K3 fibration birational to our $X$ in
family~75 (see~\S\ref{sec:constr}); proving a technical result,
Theorem~\ref{thm:75aux}, using exclusion arguments (\S\ref{sec:excl}); and
deriving Theorem~\ref{thm:75main} from Theorem~\ref{thm:75aux}
(in~\S\ref{sec:class}) using an analogue of the
Noether--Fano--Iskovskikh inequalities (\ref{nfi}) together with an
adaptation of the framework of~\cite{Ch00}.  We now make some comments on
each of these three.

\begin{emp}
In \S\ref{sec:constr} we show that the projection $(x^4,y)
\colon X \dashto \PP^1$ is indeed a K3 fibration, after resolution of
indeterminacy.  The construction is Mori-theoretic: we make a Kawamata
blowup of $X$ and play out the two ray game.  We also outline constructions
of the elliptic fibrations in
Theorems~\ref{thm:34main}(b),~\ref{thm:88main}(b) and~\ref{thm:90main}(b).
\end{emp}

\begin{comments_on_excl_sec} \label{comments_on_excl_sec}
First let us state the technical theorem mentioned above,
which is proved in~\S\ref{sec:excl}.  We need the following.

\begin{notn*}
Let $X$ be a normal complex projective variety, $\HH$ a mobile linear
system on $X$ and $\alpha \in \QQ_{\ge 0}$. We denote by
$\CS(X,\alpha\HH)$ the set of centres on $X$ of valuations that are
strictly canonical or worse for $K_X + \alpha\HH$ --- that is,
\[ \CS(X,\alpha\HH) = \{\Centre_X(E) \mid a(E,X,\alpha\HH) \le 0\}. \]
Occasionally we also use $\LCS(X,\alpha\HH)$, which is defined similarly as
\[ \LCS(X,\alpha\HH) = \{\Centre_X(E) \mid a(E,X,\alpha\HH) \le
{-1}\}. \]
\end{notn*}

\begin{thm} \label{thm:75aux}
Let $X = X_{30} \subset \PP(1,4,5,6,15)_{x,y,z,t,u}$ be a general
member of family~75 of the~95.  Suppose $\HH$ is a mobile linear
system of degree $n$ on $X$ with $\KXnH$ nonterminal.  Then in fact
$\KXnH$ is strictly canonical and $\CSXnH = \{Q_1,Q_2\}$, where
$Q_1,Q_2 \sim \frac{1}{5}(1,4,1)_{x,y,t}$ are the two singularities of
$X$ on the $zu$-stratum.
\end{thm}

For an introduction to the relationship between this result and
Theorem~\ref{thm:75main}, see~\ref{comments_on_final_sec} below.  Our proof
of~\ref{thm:75aux} in~\S\ref{sec:excl} is by exclusion arguments, as
mentioned above --- e.g., we show in Theorem~\ref{thm:smpts} that no smooth
point can be a centre on $X$ of a valuation strictly canonical for $\KXnH$;
we refer to this as \emph{excluding\/} any smooth point.  As well as
\emph{absolute\/} exclusions such as this, there is an interesting
\emph{conditional\/} exclusion result for singular points,
Theorem~\ref{thm:Tmethodcond}.  Some of the exclusion results
of~\S\ref{sec:excl} are extensions of arguments in~\cite{CPR}, but
Theorem~\ref{thm:Tmethodcond} is an example of strikingly new behaviour,
and there are substantial differences in method also for exclusion of
curves~(\S\ref{subsec:curves}).  It should be noted that, though the main
aim of~\S\ref{sec:excl} is to prove Theorem~\ref{thm:75aux}, several of the
results obtained apply to many of the~95 families other than number~75 (in
particular, numbers~34,~88 and~90) --- and the techniques of proof
apply more generally still.  \cite[App.\ A]{Ry02} contains a detailed
list of results analogous to Theorem~\ref{thm:75aux} for most of
the~95 families, including many for which only a conjectural
birational classification of elliptic and K3 fibrations is currently
known.  In the case of family~34, for instance, we have the following.

\begin{thm} \label{thm:34aux}
Let $X = X_{18} \subset \PP(1,1,2,6,9)_{x_0,x_1,y,z,t}$ be a general
member of family~34 of the~95. Suppose $\HH$ is a mobile linear system
of degree $n$ on $X$ with $\KXnH$ nonterminal. Then in fact $\KXnH$ is
strictly canonical and $\CSXnH$ is either $\{C,P,Q_1,Q_2,Q_3\}$ or
$\{P\}$. Here $C = \{x_0=x_1=0\}\cap X$ is irreducible by generality
of $X$, and $P$, $Q_1$, $Q_2$ and $Q_3$ are the singularities of $X$;
$P \sim \frac{1}{3}(1,1,2)_{x_0,x_1,y}$ lies on the $zt$-stratum and
$Q_1,Q_2,Q_3 \sim \frac{1}{2}(1,1,1)_{x_0,x_1,t}$ lie on the
$yz$-stratum.
\end{thm}

\noindent Like Theorem~\ref{thm:34main}, this result is not proved in
this paper, but the techniques for doing so, and for proving
analogues for families~88 and~90, are essentially those we use
below for Theorem~\ref{thm:75aux}.
\end{comments_on_excl_sec}

\begin{comments_on_final_sec}  \label{comments_on_final_sec}
In a sense the relationship between Theorems~\ref{thm:75aux}
and~\ref{thm:75main} is obvious: the K3 fibration given by $(x^4, y)$
has the singularities $Q_1$ and $Q_2$ of~\ref{thm:75aux} as centres,
and no other set of centres is possible --- so if, say, we try to grow
an elliptic fibration birational to $X$, we must start with an
extremal extraction of $Q_1$ or $Q_2$, but then we find we have to
extract the other $Q_i$ as well.
It turns out that to make this rigorous we need some abstract machinery
--- in particular, results concerning the log Kodaira dimension of
$(X,(\frac{1}{n} + \ep)\HH)$ for small $\ep$.  It should be
noted that for many of the~95 families an argument such as the above
does not apply directly, because sets of centres do not distinguish
objects of interest: for example, \cite{Ry02}
contains many examples of $\QQ$-Fanos $X$ birational to, say, an elliptic
fibration and also a Fano with canonical singularities, and with the two
linear systems having the same $\CS$ on $X$.

To describe the approach of~\S\ref{sec:class} we consider a more abstract
setup.  Let $X$ be a Mori Fano variety (see~\ref{defns:Mfs}
below), $Z$ a variety with canonical singularities and $\Phi \colon X
\dashto Z$ a birational map. Assume
furthermore that one of the following holds.

(a) $g \colon Z \to T$ is a $K$-trivial fibration (see~\ref{defns:fibr})
with $0 < \dim T < \dim Z$.  Let $\HH_Z = g^*|H|$
be the pullback of a very ample complete linear system of Cartier
divisors on $T$ and $\HH$ its birational transform on $X$.

(b) $Z$ is a Fano variety with canonical singularities
(see~\ref{defns:Mfs}) and $\Phi$ is
not an isomorphism. Let $\HH_Z = |H|$ be a very ample complete linear
system of Cartier divisors on $Z$ and $\HH$ its birational transform on
$X$.

In either case (a) or (b), define $n \in \QQ$ by $\KXnH \qeq 0$.

\begin{nfi}[{\cite{Ch00}}] \label{nfi}
In either of the above situations, $\KXnH$ is nonterminal, that is,
strictly canonical or worse.
\end{nfi}

This is a standard result used in~\cite{Ch00} and~\cite{Is01}; but we
give a proof of it in~\S\ref{sec:class}, under the assumptions~(a),
because in this situation it follows from results we need
anyway. Clearly~\ref{nfi} is motivation enough to address
Theorem~\ref{thm:75aux} on the way towards proving
Theorem~\ref{thm:75main}, but more work is needed to complete the
proof of the latter. \S\ref{sec:class} contains the necessary
arguments; two of the propositions are results of Cheltsov (though we
give a proof of one of them), but in order to conclude we need a
rather delicate argument that traces a log Kodaira dimension through a
two ray game diagram.
\end{comments_on_final_sec}

\subsection*{What is special about families 34, 75, 88 and 90?}

It is natural to ask why we are able to prove
Theorems~\ref{thm:75main},~\ref{thm:34main},~\ref{thm:88main}
and~\ref{thm:90main} for families~75,~34,~88 and~90, but are not able to
prove similar theorems for other families out of the~95 with the same
methods.  There is no one answer to this, but the following are major
factors.
\begin{itemize}
\item General members of families~34,~75,~88 and~90 are \emph{superrigid},
i.e., there are no nonautomorphic birational selfmaps.  For nonsuperrigid
families (the majority of the~95) one obtains much bigger anticanonical
rings after Kawamata blowups of the singularities that are centres of
involutions: see~\cite{CPR}.  This makes impossible any direct
generalisation of arguments such as that in our final proof of
Theorem~\ref{thm:75main}(b).
\item As already mentioned, it frequently occurs for other families out of
the~95 that sets of centres on $X$ do not distinguish between different
birational maps to elliptic or K3 fibrations or to Fanos with canonical
singularities. \cite[Ch.\ 4]{Ry02} discusses this phenomenon in some
detail using family~5, $X_7 \subset \PP(1,1,1,2,3)$, as an extended
example.  Whilst the analogue of Theorem~\ref{thm:75main} is eventually
proved for this family in \cite{Ry02}, the proof requires complicated
exclusion arguments on blown up models of $X$, and there is no obvious way
of avoiding these.  See~\cite{CM04} for a situation that is in some ways
analogous.
\end{itemize}

The recent work of Cheltsov and Park~\cite{CP05} proves a number of
results that complement those here: they show, for example, that a
general member of any of the~95 families is birational to a K3
fibration, but do not classify the K3 fibrations so obtained; they
also prove that a general member of family~$N$ is not birational to an
elliptic fibration if and only if $N \in \{3,75,84,87,93\}$ and, for
some 23 values of $N$ not
in this set, they show that up to a birational twist of the base there
is a unique birational elliptic fibration: it is obtained by
projection $(x_0,x_1,x_2)$ onto the first three
coordinates~\cite[4.13]{CP05}.  They do not, however, prove full
analogues of our Theorem~\ref{thm:75aux} for these families; this is
one reason why they cannot classify birational K3 fibrations.  Our
methods do permit K3 fibrations to be classified, but only in cases
where whenever we blow up canonical centres not excluded by general
results we obtain small, manageable anticanonical rings.  With our
current technology this restricts us to families such as~34,~75,~88
and~90, where there are few birational maps to elliptic and K3
fibrations and no nontrivial birational maps to Fanos with canonical
singularities.

\subsection*{Conventions and assumptions}

Our notations and terminology are mostly as in, for
example,~\cite{KM}, but we list here some conventions that are
nonstandard, together with assumptions that will hold throughout.

\begin{emp}
All varieties considered are complex, and they are projective and normal
unless otherwise stated.
\end{emp}

\begin{emp}
For details on the famous~95 families, see~\cite{Fl00}
and~\cite{CPR}.  In brief, $X_d \subset \PP(1,a_1,\ldots,a_4)$ belongs
to one of the families if (1) $X$ is \emph{quasismooth}, i.e., its
singularities are all quotient singularities forced by the weighted
$\PP^4$; (2) the singularities are \emph{terminal} --- 3-fold terminal
quotient singularities are necessarily of the form
$\frac{1}{r}(1,a,r-a)$ with $r \ge 1$ and $(a,r) = 1$; and (3) $a_1 +
\cdots + a_4 = d$, so by the adjunction formula ${-K_X} = \Oh_X(1)$ is
very ample.  Whenever $X$ is a member of one of the~95 families, we
let $A = {-K_X} = \Oh_X(1)$ denote the positive generator of the class
group; moreover, if $f \colon Y \to X$ is a birational morphism then
$B$ denotes ${-K_Y}$.
\end{emp}

\begin{emp}
As in~\cite{CPR} and~\cite{Ry02}, we refer to the weighted blowup with
weights $\frac{1}{r}(1, a, r-a)$ of a 3-fold terminal quotient singularity
$\frac{1}{r}(1, a, r-a)$ with $r \ge 2$ and $(a, r) = 1$ as the
\emph{Kawamata blowup\/} --- see~\cite{Ka96}, the main theorem of which is
reproduced here as Theorem~\ref{thm:Ka}.
\end{emp}

\begin{defns} \label{defns:Mfs}
A \emph{Mori fibre space\/} $f \colon X \to S$ is a Mori extremal contraction
of fibre type, that is, $\dim S < \dim X$. This means that $X$ and $S$
are projective varieties, $X$
has $\QQ$-factorial, terminal singularities, $f_*\Oh_X = \Oh_S$,
$\rho(X/S) = 1$ and ${-K_X}$ is $f$-ample. If $S = \{*\}$ is a point
then $X$ is a \emph{Mori Fano variety}. We use the term \emph{Fano
variety\/} more generally to refer to any normal, projective variety $X$ with
$-K_X$ ample and $\rho(X) = 1$.
\end{defns}

\begin{defns} \label{defns:fibr}
Let $Z$ be a normal projective variety with canonical singularities.
A \emph{fibration\/} is a morphism $g \colon Z \to T$ to another normal
projective variety $T$ such that $\dim T < \dim Z$ and $g_*\Oh_Z =
\Oh_T$. We say the fibration is \emph{$K$-trivial\/} if and only if $K_Z C =
0$ for every contracted curve~$C$. $g$ is an \emph{elliptic
  fibration}, resp.\ a \emph{K3 fibration}, if and only if its general
fibre is an elliptic curve, resp.\ a K3 surface.
\end{defns}

\begin{emp}
Usually when we write an equation explicitly or semi-explicitly in
terms of coordinates we omit scalar coefficients of monomials; this is the
`coefficient convention'.
\end{emp}

\subsection*{Acknowledgments}

Most of the techniques in this paper, and some of the
theorems, are from my PhD thesis,~\cite{Ry02}.  I would like to thank
my supervisor, Miles Reid, and also Gavin Brown, Alessio Corti and
Hiromichi Takagi, for their help and generosity with their ideas; I
would also like to thank Ivan Cheltsov for his helpful comments during
the final stages of preparing this paper.  My PhD studies were
supported financially by the British EPSRC.

\section{Constructions} \label{sec:constr}

\subsection{K3 fibrations}

The following observation does not apply to family~75, our main object
of study, but it does apply to families~34 and~88; in any case, it needs to
be noted because it describes all `easy' K3 fibrations birational to
members of the famous 95 --- cf.~Lemma~\ref{lem:75tworay} for the `hard'
case.

\begin{prop} \label{prop:easyK3s}
Let $X = X_d \subset \PP(1,1,a_2,a_3,a_4)$ be general in one of the
families with $a_1 = 1$ and $a_2 > 1$. Then a general fibre $S$ of
$\pi = (x_0,x_1) \colon X \dashto \PP^1$ is a quasismooth Du Val K3
surface and, setting $\PPP$ to be the pencil $\left<x_0,x_1 \right>$,
we have
\[ \CS(X,\PPP) = \{C,P_1,\ldots,P_r\}, \]
where $C$ is the curve $\{x_0 = x_1 = 0\} \cap X$, which is
irreducible by generality of $X$, and $P_1,\ldots,P_r$ are all the
singularities of $X$.
\end{prop}

\begin{proof}[\textsc{Proof}]
Because $S$ is a general element of $|\Oh_X(1)|$, it is certainly
quasismooth. The adjunction formula for $S_d \subset
\PP(1,a_2,a_3,a_4) =: \PP$ gives $K_S = 0$ and the cohomology long
exact sequence from
\[ 0 \to \II_{S,\PP}=\Oh_{\PP}(-d) \to \Oh_{\PP} \to \Oh_S \to 0, \]
with the standard cohomology results for weighted projective space,
gives $h^1(S,\Oh_S) = 0$. Therefore $S$ is a quasismooth Du Val K3
surface.

Let $f \colon Y \to X$ be the blowup of the ideal sheaf $\II_{C,X}$ of
$C$ and $E \subset Y$ the unique exceptional divisor of $f$ which
dominates $C$. Then clearly $m_E(\PPP) = a_E(K_X) = 1$, so $C \in
\CS(X,\PPP)$. The fact that $P_1,\ldots,P_r \in \CS(X,\PPP)$ is a
consequence of Corollary~\ref{cor:Kalemma:2} of Kawamata's Lemma
(below). Therefore
\[ \CS(X,\PPP) \supset \{C,P_1,\ldots,P_r\}; \]
the reverse inclusion follows from Theorem~\ref{thm:smpts}.
\end{proof}

\begin{rk}
If $X = X_d \subset \PP(1,1,1,a_3,a_4)$ is general in a family with $a_1 =
a_2 = 1$ then clearly a general element $S$ of any pencil $\PPP \subset
|\Oh_X(1)|$ is a Du Val K3 surface, provided one can prove it is
quasismooth.  This can be fiddly, the problem being that while $X$ is
general and $S \in \PPP$ is general, $\PPP$ must be able to be
\emph{any\/} pencil inside $|\Oh_X(1)|$.
\end{rk}

In contrast to the situation considered above, for families $X_d \subset
\PP(1,a_1,\ldots,a_4)$ with $a_1 > 1$ it is not immediately clear whether
there exist K3 fibrations birational to $X$: to construct them, or at least
to make sense of the construction, we need Mori theory.  Here we consider
only family~75, but the technique applies to many other families ---
see~\cite{Ry02} --- and, in particular, to family~90, the subject of our
Theorem~\ref{thm:90main}.

\begin{lemma} \label{lem:75tworay}
Let $X=X_{30} \subset \PP(1,4,5,6,15)_{x,y,z,t,u}$ be a general
member of family~75 of the~95; for the notations $P,Q_1,Q_2,R_1,R_2$ for
the singularities of $X$, see Theorem~\ref{thm:75aux}. Let $Q_i \in X$
be either $Q_1$ or $Q_2$ and $f \colon Y \to X$ the Kawamata blowup
of $Q_i$.
\begin{itemize}
\item[(1)] Let $R \subset \NEbar Y$ be the ray with $\cont_R = f$. Then the
other ray $Q \subset \NEbar Y$ is contractible and its contraction $g
= \cont_Q \colon Y \to Z$ is antiflipping.
\item[(2)] The antiflip $Y \dashrightarrow Y'$ of $g$ exists and $Y'$ has
canonical singularities.
\item[(3)] Let $Q' \subset \NEbar Y'$ be the ray whose contraction is $g'
\colon Y' \to Z$ and $R' \subset \NEbar Y'$ the other ray. Then $R'$
is contractible and its contraction $f' \colon Y' \to \PP^1$, which is in fact the
anticanonical morphism $\fie_{|{-4K_{Y'}}|} = \fie_{|4B'|}$, is a K3
fibration --- that is, a general fibre $T'$ of $f'$ has Du Val
singularities, $K_{T'} =0$ and $h^1(T',\Oh_{T'}) = 0$.
\end{itemize}
It follows that the total composite $X \dashto \PP^1$ of the two ray game we
have played, illustrated below, is $\pi = (x,y) \colon X \dashto
\PP(1,4) = \PP^1$.
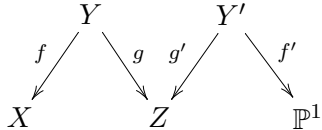
\[ \xymatrix@C=0.4cm{
 & Y \ar[dl]_f \ar[dr]^g & & Y' \ar[dl]_{g'} \ar[dr]^{f'} & \\
 X & & Z & & \PP^1 \\
} \]
Therefore $R(Y,B) = R(Y',B') = k[x,y]$.
\end{lemma}

\begin{proof}[\textsc{Proof of \ref{lem:75tworay}}]
(1) The first part of the following argument --- showing that the
curve $C$ defined below generates $Q \subset \NEbar Y$ --- is one case
of \cite[5.4.3]{CPR}; but in \cite{CPR} this point $Q_i$ is excluded
as a maximal centre by the test class method, so there is no need for
the two ray game to be played out.

By generality of $X$, the curve $\{x=y=0\} \cap X$ is
irreducible. $f \colon Y \to X$ is the $\frac{1}{5}(1,4,1)_{x,y,t}$
weighted blowup of the $\frac{1}{5}(1,4,1)_{x,y,t}$ point
$Q_i$. Let $S \in \left|B\right|$ be the unique effective surface and
$T \in \left|4B\right|$ a general element, where as always $B =
{-K_Y}$.  One can check explicitly, by looking at the three affine
pieces of $Y$ locally over $Q_i$, that $C := S \cap T$ is an
irreducible curve inside $Y$; this uses the generality of
$X$. We also know that
\[ B^3 = A^3 - \frac{1}{ra(r-a)} = \frac{1}{60} - \frac{1}{20} < 0, \]
so in particular $BC = 4B^3 < 0$. It follows that the ray $Q \subset
\NEbar Y$ is generated by $C$ --- indeed, suppose this is not the
case; then $C$ is in the interior of the 2-dimensional cone $\NEbar
Y$, so we can pick an effective 1-cycle $\sum_{i = 0}^p \alpha_i C_i$ that lies
strictly between $Q$ and the half-line generated by $C$. This 1-cycle
is $B$-negative, because $R$ is $B$-positive (i.e., $K$-negative) and
$BC < 0$; but $\Bs\left|4B\right|$ is supported on $C$ and therefore
one of the $C_i$, say $C_0$, is in fact $C$. The geometry of the cone
now implies that, after we subtract off $\alpha_0 C_0$, $\sum_{i =
  1}^p \alpha_i C_i$ is again strictly between $Q$ and the half-line
generated by $C$; so we can repeat the argument to deduce that some
other $C_i$ is $C$ --- and of course this contradicts our initial (implicit)
assumption that the $C_i$ were distinct. This argument has also shown
that $C$ is the only irreducible curve in the ray $Q$.

We show $Q$ is contractible with a general Mori-theoretic trick ---
afterwards, clearly, $g = \cont_Q$ is antiflipping, because $C$ is the
only contracted curve and $K_Y C = {-BC} > 0$. Firstly note that $S$
has canonical singularities and so in particular is klt. We apply
Shokurov's inversion of adjunction (see \cite[5.50]{KM}) to deduce
that the pair $(Y,S)$ is plt. (In fact $K_Y + S$ is Cartier so the
log discrepancy of any valuation for $(Y,S)$ is
an integer, and therefore $(Y,S)$ is canonical.) But plt is an open
condition so for $\ep \in \QQ_{>0}$, $\ep \ll 1$, the pair
$(Y,S + \ep T)$ is plt as well. Now
\[ (K_Y + S + \ep T)C = \ep TC = \ep B^3 < 0,  \]
so $Q$ is contractible by the (log) contraction theorem for $(Y,S +
\ep T)$.

(2) This follows immediately from Mori's result
\cite[20.11]{FA}: $S$ and $T$ are effective divisors, $T \sim 4S$ and
$T \cap S = C$ is precisely the exceptional set of $g$. Since
$S \sim B = {-K_Y}$ the antiflip of $g$ is precisely its `opposite
with respect to $S$', to use the language of \cite{FA}. This can be
constructed as the normalisation of the closure of the image of
\[ g \times \pi_Y \colon Y \dashrightarrow Z \times \PP^1, \]
where $\pi_Y = \pi \circ f \colon Y \dashrightarrow \PP(1,4) = \PP^1$
corresponds to the pencil $\left<4S,T\right> = \left|4B\right|$.

Alternatively, since we have already observed that the ray $Q$ is
$(K_Y + S + \ep T)$-negative, (2) follows from Shokurov's general
result that log flips of lc pairs exist in dimension 3.

(3) From the construction of the antiflip, the transforms $S'$ and
$T'$ of $S$ and $T$ on $Y'$ are disjoint, so
$\Bs\left|{-4K_{Y'}}\right| = \Bs\left|4B'\right| =
\emptyset$. Therefore $f' = \cont_{R'}$ exists and is (the Stein
factorisation of) $\fie_{\left|4B'\right|}$; $T'$ is a general fibre
of $f'$. In the diagram below, $\nu \colon U \to g(T)$ is the
normalisation of $g(T)$.
\[ \xymatrix@C=0.1cm{
 & Y & \supset & T \ar[dl]_f \ar@/_0.3cm/[ddrr]_g \ar[drr]^h  & & & &
 T'  \ar[dll]_{h'} \ar@/^0.3cm/[ddll]^{g'} \ar[dr]^{f'}  & \subset &
 Y'& \\  X & \supset & T_X & & & U  \ar[d]^{\nu} & & & \{*\} & \subset
 & \PP^1 \\  & & & & & g(T) & & & & & \\
} \]
Now $K_T = (K_Y + T)|_T = 3B|_T = 3C$, so $K_U = h_*(3C) =
0$. Furthermore, the Leray spectral sequence for $f\colon T\to T_X$,
\[ 0 \to 0 = H^1\left(T_X,\Oh_{T_X}\right) \to
 H^1\left(T,\Oh_T\right) \to H^0\left(T_X,R^1f_*\Oh_T\right) = 0 , \]
shows that $h^1(T,\Oh_T) = 0$ (here $h^0(T_X,R^1f_*\Oh_T) = 0$
because the singularity $\frac{1}{5}(1,1)$ of $T_X$ at $Q_i$ is
rational); and now the Leray spectral sequence for $h \colon T \to U$,
\[ 0 \to H^1\left(U,\Oh_U\right) \to H^1\left(T,\Oh_T\right) = 0 , \]
shows that $h^1(U,\Oh_U) = 0$.
The only thing left to do is to show that $U$ has Du Val singularities
--- it is then clear that $T'$ is a Du Val K3 surface because
\[ K_{T'} = \left(K_{Y'} + T'\right)|_{T'} = 3B'|_{T'} = 0 , \]
so the minimal resolution $\widetilde{U} \to U$ of $U$ factors through
$h' \colon T' \to U$.

To show $U$ has Du Val singularities we observe first that $T$ has two
singularities: a $\frac{1}{5}(1,1)_{x,t}$ point over $Q_j \in T_X$,
where $\{i,j\} = \{1,2\}$, and a $\frac{1}{3}(1,2)_{x,z}$ point over
$P \in T_X$.  The curve $C$ passes through both and is locally defined
by $x = 0$ in a neighbourhood of each.  Let $\xi \colon \widetilde{T}
\to T$ be the minimal resolution of $T$, with exceptional curves
$E_1$, $E_2$ lying over $P$ and $E_3$ lying over $Q_j$, where $E_1^2 =
E_2^2 = {-2}$, $E_1 E_2 = 1$, $E_1 \widetilde{C} = 0$, $E_2
\widetilde{C} = 1$, $E_3^2 = {-5}$ and $E_3 \widetilde{C} = 1$ (here
$\widetilde{C}$ is the birational transform of $C$ on
$\widetilde{T}$).  We can calculate that $(\widetilde{C})^2 < 0$ and
$K_{\widetilde{T}}\widetilde{C} < 0$, from which it follows that
$\widetilde{C}$ is a $-1$-curve.  (To do this, first show that $C^2_T
< 0$ and $K_T C < 0$; then express $\widetilde{C}$ as the pullback
$\xi^*C$ minus nonnegative multiples of $E_1,E_2,E_3$; then use the
fact that $(E_i E_j)_{i,j = 1}^2$ is negative definite,
$K_{\widetilde{T}}(\xi^*C) = K_T C < 0$ and $K_{\widetilde{T}}$ is
$\xi$-nef by minimality of $\xi$.)

Now the minimal resolution $\widetilde{U}$ of $U$ is obtained from
$\widetilde{T}$ by running a minimal model program over $U$ --- so we
start by contracting $\widetilde{C}$, which is exceptional over $U$,
and then $E_2$, which has become a $-1$-curve, and finally $E_1$.
$E_3$ is left as a $-2$-curve.  This MMP can be summarised by
\[ (2,2,1,5) \to (2,1,4) \to (1,3) \to (2). \]
It follows that $U$ has a single $A_1$ Du Val singularity, as
required.
\end{proof}

\subsection{Elliptic fibrations}

Our main aim is to prove Theorem~\ref{thm:75main}, part of which is the
statement that there is no elliptic fibration birational to a general $X$
in family~75.  In this subsection, however, we digress briefly to discuss
the construction of elliptic fibrations birational to hypersurfaces in
other families out of the~95.  In particular we are concerned with
families~34,~88 and~90, the subjects of
Theorems~\ref{thm:34main},~\ref{thm:88main} and~\ref{thm:90main}.

\begin{ex}
Let $X = X_{18} \subset \PP(1,1,2,6,9)_{x_0, x_1, y, z, t}$ be a general
member of family~34 of the~95.  We claim that
\begin{itemize}
\item[(a)]  the projection $\pi = (x_0, x_1, y) \colon X \dashto
\PP(1,1,2)$ gives an elliptic fibration birational to $X$ after resolution
of indeterminacy; and
\item[(b)]  the indeterminacy may be resolved as shown below, where $f
\colon Y \to X$ is the Kawamata blowup of the unique singularity $P \sim
\frac{1}{3}(1,1,2)_{x_0,x_1,y}$ on the $zt$-stratum, and $\pi_Y := \pi
\circ f$ is the anticanonical morphism $\fie_{|2B|}$, $B = {-K_Y}$.
\[ \xymatrix@C=1.6cm{
Y \ar[rd]^(.4){\pi_Y} \ar[d]_{f} & \\
X \ar@{-->}[r]^(.4){\pi} & \PP(1,1,2) \\
} \]
\end{itemize}

\begin{proof}[\textsc{Proof}]
(a) This is easy to see: a general fibre of the rational map $\pi$ is a
curve $E_{18} \subset \PP(1,6,9)_{x,z,t}$ and, when we write down the
Newton polygon of the defining equation of such a curve, there is a unique
internal monomial (namely $x^3 z t$, the vertices being $t^2$,
$x^{18}$ and~$z^3$).  Consequently $E_{18}$ is birational to an
elliptic curve, by standard toric geometry.

(b) The linear system $\LL$ defining $\pi$ is $|2A| = \left< x_0^2, x_1^2,
y \right>$ and one can calculate directly that its birational transform
$\LL_Y$ on $Y$ is free.  Furthermore
\[ \LL_Y = f^*\LL - \textstyle\frac{2}{3} E \sim 2B \]
and $\LL_Y$ is clearly a complete linear system.
\end{proof}
\end{ex}

As well as exhibiting one elliptic fibration birational to a general
hypersurface in family~34 (actually, according to
Theorem~\ref{thm:34main}, the only such), this example shows us one
way to look for elliptic fibrations when we consider other families:
namely, find a singular point $P$ with $B^3 = 0$ (for $B = {-K_Y}$, $f
\colon Y \to X$ the Kawamata blowup of $P$), take the anticanonical
morphism on $Y$ (if $B$ is eventually free) and see if it maps to a
surface.  It turns out that this method works for families~88 and~90,
as well as~34, so as far as the present paper is concerned, we are
done.  There are, however, other ways in which elliptic fibrations
occur birational to hypersurfaces in the~95 families.  Here is a brief
list; see~\cite{Ry02} for more details.
\begin{itemize}
\item  Sometimes elliptic fibrations have more than one singular point
of $X$ in $\CSXnH$; they can be constructed by blowing up all these
points and taking the anticanonical morphism.
\item  {\sloppy It can also occur that $\CSXnH$ consists of only one
singular point $P$, but \linebreak $\CS\left( Y, \frac{1}{n}\HH_Y
\right) \ne \emptyset$, where $Y = \B_P X$ is the Kawamata blowup; in
this case further blowups of $Y$ are necessary before taking the
anticanonical morphism.}
\item  Finally, there are examples of elliptic fibrations with a curve
in $\CSXnH$ --- but only for families~1 and~2; see~\cite{Ry02}.
\end{itemize}

\section{Exclusions, absolute and conditional} \label{sec:excl}

For an initial introduction to the contents of this section,
see~\ref{comments_on_excl_sec}.  We divide the material for proving
Theorem~\ref{thm:75aux} into three subsections:
in~\S\ref{subsec:curves} we show that all curves are excluded
absolutely and in~\S\ref{subsec:smpts} we prove the corresponding
result for smooth points; finally in~\S\ref{subsec:singpts} we deal
with singular points.  Much of the material in this section applies
more widely than to family~75: for example, out of all the singular
points on members of the~95 families, we show that those satisfying a
certain condition turn out to be excluded absolutely, while those
satisfying a different (but closely related) condition are excluded
\emph{conditionally} --- that is, we prove that if they belong to
$\CSXnH$ for some $\HH$ then other centres must also exist in
$\CSXnH$.

\subsection{Curves} \label{subsec:curves}

First we state the main curve exclusion theorem proved in \cite{Ry02}.

\begin{thm}[{\cite[Curves Theorem A]{Ry02}}] \label{thmA}
Let $X = X_d \subset \PP = \PP(1,a_1,a_2,a_3,a_4)$ be a general
hypersurface in one of the 95 families and $C \subset X$ a reduced,
irreducible curve. Suppose $\HH$ is a mobile linear system of degree
$n$ on $X$ such that $K_X + \frac{1}{n}\HH$ is strictly canonical and
$C \in \CSXnH$. Then there exist two linearly independent forms
$\ell,\ell'$ of degree $1$ in $(x_0,\ldots,x_4)$ such that
\begin{equation} C \subset \{\ell = \ell' = 0\} \cap X. \label{CinPi}
\end{equation} \end{thm}

We do not reproduce the proof in full here, but instead restrict
ourselves to the case we need, namely $X_{30} \subset
\PP(1,4,5,6,15)$.  This follows immediately from the following lemmas,
the first of which is standard.

\begin{lemma} \label{lem:degC}
Let $X$ be any hypersurface in one of the 95 families and $C \subset
X$ a curve, reduced but possibly reducible. Suppose $\HH$ is a mobile
linear system of degree $n$ on $X$ such that $K_X + \frac{1}{n}\HH$ is
strictly canonical and each irreducible component of $C$ belongs to
$\CSXnH$. Then $\deg C = AC \le A^3$.
\end{lemma}

\begin{lemma} \label{curvesinwps1}
Let $X = X_d \subset \PP = \PP(1,a_1,a_2,a_3,a_4)$ be a hypersurface
in one of the families with $a_1 > 1$ and suppose that either
\begin{itemize}
\item[(a)] $d < a_1a_4$ or
\item[(b)] $d < a_2a_4$ and the curve $\{x = y = 0 \}\cap
X$ is irreducible (which holds for general $X$ in a family with $a_1 >
1$ by Bertini's theorem). 
\end{itemize}
Then any curve $C \subset X$ that is not
contracted by $\pi_4 \colon X \dashrightarrow \PP(1,a_1,a_2,a_3)$
has $\deg C > A^3$. Consequently $C$
is excluded absolutely by Lemma \ref{lem:degC}.
\end{lemma}

For the proofs of Lemmas~\ref{lem:degC} and~\ref{curvesinwps1}, see
below.  It is straightforward to check that they imply the following.

\begin{cor} \label{cor:75curves}
Let $X = X_{30} \subset \PP(1,4,5,6,15)_{x,y,z,t,u}$ be any
(quasismooth) member of family~75. Then Theorem~\ref{thmA} holds for
$X$, that is, no reduced, irreducible curve $C \subset X$ can belong
to $\CSXnH$ for any mobile system $\HH$ of degree $n$ with $\KXnH$
strictly canonical.
\end{cor}

\begin{proof}[\textsc{Proof of Lemma \ref{lem:degC}}]
Let $s$ be a natural number such that $sA$ is
Cartier and very ample, and pick general members $H,H' \in
\HH$. Now by assumption
\[ \mult_{C_i}(H) = \mult_{C_i}(H') = n \]
for each irreducible component $C_i$ of $C$, so for a general $S \in
\left|sA \right|$
\[ A^3 sn^2 = SHH' \ge sn^2 AC = sn^2 \deg C, \]
which proves $\deg C \le A^3$.
\end{proof}

\begin{proof}[\textsc{Proof of Lemma \ref{curvesinwps1}}]
We prove this lemma under the additional assumption that $(a_1,a_2) =
1$ --- which is the case for family 75 ($a_1 = 4$, $a_2 = 5$); if
$(a_1,a_2) > 1$ then a little trick, described in \cite{Ry02}, is
needed.

Suppose that $C \subset X$ has $\deg C \le A^3$ and is not contracted
by $\pi_4$; let $C' \subset \PP(1,a_1,a_2,a_3)$ be the
set-theoretic image $\pi_4(C)$. Note that $\deg C' \le \deg C$ ---
indeed, if $H$ denotes the hyperplane section of
$\PP(1,a_1,\ldots,a_4)$ and $H'$ that of $\PP(1,a_1,a_2,a_3)$, we pick
$s \ge 1$ such that $\left|sH\right|$ and $\left|sH'\right|$ are very
ample, and calculate that
\begin{eqnarray*}
s\deg C & = & (sH)C = \pi_4^*(sH')C \\
        & = & sH'(\pi_4)_*C = srH'C' = sr\deg C' \ge s\deg C',
\end{eqnarray*}
where $r \ge 1$ is the degree of the induced morphism $\pi_4|_C \colon
C \to C'$. So in fact $\deg C$ is a multiple of $\deg C'$ by the
positive integer $r$.
(The point of $\left|sH\right|$ being very ample is that we can move it
away from $P_4$, where $\pi_4$ is undefined, and apply the projection
formula to the morphism $\pi_4|_{\PP(1,a_1,\ldots,a_4) \smallsetminus
\{P_4\}}$.)

Now form the diagram below.
\[ \xymatrix@C=0.1cm{
 C \ar[d]    & \subset & \PP(1,a_1,a_2,a_3,a_4) \ar@{-->}[d]^{\pi_4}\\
 C' \ar[d]   & \subset & \PP(1,a_1,a_2,a_3) \ar@{-->}[d]^{\pi_3}\\
 \{*\}       & \subset & \PP(1,a_1,a_2)\\
} \]
$C'$ is contracted by $\pi_3$
 --- indeed, if its image were a curve $C''$ we would have
\[ \deg C'' \le \deg C' \le \deg C \le A^3,  \]
but $A^3 = d/(a_1a_2a_3a_4) < 1/(a_1a_2)$,
since $d < a_3a_4$ in either case (a) or (b), and on the other hand
$1/(a_1a_2) \le \deg C''$
simply because $C'' \subset \PP(1,a_1,a_2)$ --- contradiction.

Now by our extra assumption $(a_1,a_2) = 1$, the point $\{*\} \subset
\PP(1,a_1,a_2)$ is, up to coordinate change, one of
\[ \{y = z = 0\}, \quad \{y^{a_2} + z^{a_1} = x = 0\}, \quad \{x = z =
0\} \quad \mbox{and} \quad \{x = y = 0\}, \]
using the coefficient convention in $y^{a_2} + z^{a_1} = 0$. It
follows that the curve $C' \subset \PP(1,a_1,a_2,a_3)$ is defined by
the same equations. In the first case, this means that
$\deg C' = 1/a_3 > d/(a_1a_2a_3a_4) = A^3$,
contradiction. In the second case $\deg
C' = 1/a_3$ again, because $C' \iso \{y^{a_2} + z^{a_1} = 0\}
\subset \PP(a_1,a_2,a_3)$ passes only through the singularity
$(0,0,1)$, using $(a_1,a_2) = 1$ --- so we obtain a contradiction as in
the first case. In the case $C' = \{x = z = 0\}$, we have $\deg C' =
1/(a_1a_3)$ and we easily obtain a contradiction from $a_2a_4 >
d$. In the final case, $C' = \{x = y = 0\}$, if the assumptions in
part (a) of the statement hold then we have
\[ \deg C' = 1/(a_2a_3) > d/(a_1a_2a_3a_4) = A^3, \]
contradiction; while if the assumptions in part (b) hold then
\[ C = \{x = y = 0 \}\cap X \]
(because the right hand side is irreducible), but
\[ \deg\left(\{x = y = 0 \}\cap X \right) = a_1A^3 > A^3,  \]
since we also assumed $a_1 > 1$ --- contradiction.
\end{proof}

Note that for the case $a_1 = 1, a_2 > 1$ there is the following
analogue of Lemma~\ref{curvesinwps1} --- see~\cite[3.3]{Ry02} for the
proof, which is similar to but much shorter than the one we have just
seen.

\begin{lemma} \label{curvesinwps2}
Let $X = X_d \subset \PP = \PP(1,1,a_2,a_3,a_4)$ be a hypersurface in
one of the families with $a_1 = 1$ and $a_2 > 1$; suppose that $d <
a_2a_4$. Then any curve $C \subset X$ that is not contracted by
$\pi_4$ and that satisfies $\deg C \le A^3$ is contained in $\{x_0 =
x_1 = 0\} \cap X$.
\end{lemma}

\noindent When $a_1 = a_2 = 1$, however, the situation is different
and more work is required to prove sufficiently strong results.

We note also for completeness that in the case of family~75
we have $P_4 = P_u \not\in X$ so the question of whether $C$ is
contracted by $\pi_4 \colon X \dashto \PP(1,a_1,a_2,a_3)$ never
arises.  For families that do contain curves contracted by $\pi_4$,
\cite[3.5]{Ry02} shows that in almost all cases these curves are of
degree greater than $A^3$, so they are excluded by
Lemma~\ref{lem:degC}.  The remaining few families to which this result
does not apply are also dealt with in \cite{Ry02}.

\begin{crvs_fams34etc}
We have dealt with curves inside a general $X$ in family~75.  The
arguments presented above, or variants of them, deal also with curves
in general members of families~88 and~90, given the following observation:
family~88 has $a_1 = 1$ and $a_2 > 1$, and consequently there is a
birational K3 fibration obtained by the projection $(x_0, x_1) \colon
X \dashto \PP^1$.  The corresponding pencil $\PPP = \left< x_0, x_1
\right>$ has $C \in \CS(X, \PPP)$, where $C$ is the curve $\{ x_0 =
x_1 = 0 \} \cap X$.  The important point is that, by taking $X$
general in its family, $C$ is \emph{irreducible}.  If it were not, and
we had
\[ \{ x_0 = x_1 = 0 \} \cap X = C_0 \cup \cdots \cup C_r, \]
we would have to prove a conditional exclusion result of the form `if
$C_0 \in \CSXnH$ then all $C_i \in \CSXnH$'.  This can be done ---
see~\cite[3.12]{Ry02} --- but we avoid it using our generality
assumption on $X$.

The case of family~34 is more problematic because
Lemma~\ref{curvesinwps2} above does not apply.  We omit the argument
for curve exclusion for this family; it can be found
in~\cite[\S3.4]{Ry02}.
\end{crvs_fams34etc}

\subsection{Smooth points} \label{subsec:smpts}

In this subsection we present a proof of the following theorem.

\begin{smptsthm} \label{thm:smpts}
Let $X = X_d \subset \PP = \PP(1,a_1,\ldots,a_4)$ be a general
hypersurface in one of families $3,4,\ldots,95$ and $P \in X$ a smooth
point. For any $n \in \ZZ_{\ge1}$ and any mobile linear system $\HH$
of degree $n$ on $X$ we have $P \not\in \CSXnH$.
\end{smptsthm}

The proof closely follows the argument used in \cite{CPR} to exclude
smooth points as maximal centres of a pair $\XnH$. First we need to
quote some theoretical results.

\begin{thm}[Shokurov's inversion of adjunction] \label{invadj}
Let $P \in X$ be the germ of a smooth 3-fold and $\HH$ a mobile linear
system on $X$. Assume  $n \in \ZZ_{\ge 1}$ is such that $P \in
\CS\left(X,\frac{1}{n}\HH\right)$. Then for any normal irreducible divisor $S$
containing $P$ such that $\HH|_S$ is mobile we have $P \in
\LCS\left(S,\frac{1}{n}\HH|_S\right)$.
\end{thm}
For a readable account of the proof see \cite[5.50]{KM}.  Note that
under the given assumptions, but
without assuming $\left.\HH\right|_S$ is mobile, \cite[5.50]{KM} says
$K_S + \frac{1}{n}\HH|_S$ is not klt
near $P$, which does not preclude the centre on $S$ of the relevant
valuation being a \emph{curve\/} containing $P$ rather than $P$
itself. This curve would of course be in $\Bs(\HH)$, so the problem is
eliminated by assuming $\HH|_S$ is mobile --- and it will be
clear we may assume this in our application. With the extra assumption
$K_S + \frac{1}{n}\HH|_S$ is not plt near $P$, as required.

\begin{thm}[Corti] \label{Cothm}
Let $\left(P \in \Delta_1 + \Delta_2 \subset S  \right) \iso \left(0
\in \{xy = 0 \} \subset \CC^2  \right)$  be the analytic germ of a
normal crossing curve on a smooth surface; let $\LL$ be a mobile
linear system on $S$ and $L_1,L_2 \in \LL$ general members. Suppose
there exist $n \in \ZZ_{\ge 1}$ and $a_1,a_2 \in \QQ_{\ge 0}$ such
that
\[ P \in \LCS\left(S, (1-a_1)\Delta_1 + (1-a_2)\Delta_2 +
\textstyle\frac{1}{n}\LL \right).  \]
Then
\[ (L_1 \cdot L_2)_P \ge \left\{
\begin{array}{l@{\quad\mbox{if}\hspace{1.5ex}}l} 4a_1a_2n^2 & a_1 \le
1 \mbox{\hspace{1ex} or \hspace{1ex}} a_2 \le 1; \\ 4(a_1+a_2-1)n^2
& \mbox{both\hspace{1.5ex}} a_1,a_2 > 1.  \end{array} \right. \]
\end{thm}
This is proved as in the original \cite{Co00}, but replacing `log
canonical' by `purely log terminal' and strict inequalities by $\le$
or $\ge$ as appropriate. Now combining Theorems \ref{invadj} and
\ref{Cothm} we obtain the following.

\begin{cor} \label{multthm}
Let $P \in X$ be the germ of a smooth 3-fold and $\HH$ a mobile linear
system on $X$. Assuming as in Theorem \ref{invadj} that $P \in
\CSXnH$ for some $n \in \ZZ_{\ge 1}$ we have
\[ \mult_P(H \cdot H') \ge 4n^2 \]
where $H,H' \in \HH$ are general and $H \cdot H'$ is their
intersection cycle.
\end{cor}

Now we need to borrow an additional result from \cite{CPR}. First we
recall the following definition.

\begin{defn}[cf. {\cite[5.2.4]{CPR}}] \label{defn:Gammaisol}
Let $L$ be a Weil divisor class in a 3-fold $X$ and $\Gamma \subset X$
an irreducible curve or a closed point. We say that \emph{$L$ isolates
$\Gamma$}, or is a \emph{$\Gamma$-isolating class\/}, if and only if
there exists $s \in \ZZ_{\ge1}$ such that the linear system
$\LL_{\Gamma}^s := \left|\II_{\Gamma}^s(sL) \right|$ satisfies
\begin{itemize}
\item $\Gamma \in \Bs \LL_{\Gamma}^s$ is an isolated component (i.e.,
in some neighbourhood of $\Gamma$ the subscheme $\Bs \LL_{\Gamma}^s$
is supported on $\Gamma$); and
\item if $\Gamma$ is a curve, the generic point of $\Gamma$ appears
with multiplicity 1 in $\Bs \LL_{\Gamma}^s$.
\end{itemize}
\end{defn}

\begin{thm}[{\cite[5.3.1]{CPR}}] \label{Pisol}
Let $X = X_d \subset \PP(1,a_1,\ldots,a_4)$ be a general hypersurface
in one of families $3,4,\ldots,95$ and $P \in X$ a smooth point. Then
for some positive integer $l < 4/A^3$ the class $lA$ is
$P$-isolating.
\end{thm}

\begin{rk} \label{Pisolrk}
\cite[5.3.1]{CPR} says $l \le 4/A^3$, but the
statement there is for all the families except number 2 --- that is,
including number 1 --- and a trivial check shows that in fact $l =
4/A^3$ only for number 1.
\end{rk}

\begin{proof}[\textsc{Proof of Theorem \ref{thm:smpts}}]
We know that $P \in \Bs\left|\II_P^s(slA) \right|$ is an isolated
component for some $l,s \in \ZZ_{\ge1}$ with $l < 4/A^3$. Take
a general surface $S \in \left|\II_P^s(slA) \right|$ and general
elements $H,H' \in \HH$. If we assume that $P$ belongs to $\CSXnH$
then Corollary \ref{multthm} tells us that
$\mult_P \left(H\cdot H' \right) \ge 4n^2$, so
\[ {S \cdot H \cdot H'} \ge \left(S \cdot H \cdot H' \right)_P \ge
4sn^2.  \]
But we know
\[ {S \cdot H \cdot H'} = sln^2A^3 < \frac{4}{A^3}sn^2A^3 = 4sn^2,  \]
contradiction.
\end{proof}

\subsection{Singular points} \label{subsec:singpts}

The following two results are fundamental.

\begin{thm}[Kawamata, \cite{Ka96}] \label{thm:Ka}
Let $P \in X \iso \frac{1}{r}(1,a,r-a)$, with  $r \ge 2$ and $(a,r) = 1$,
be the germ of a 3-fold terminal quotient singularity, and
\[ f \colon (E \subset Y) \to (\Gamma \subset X)  \]
a divisorial contraction such that $P \in \Gamma$ (so $Y$ has terminal
singularities, $\Exc f = E$ is an irreducible divisor and $-K_Y$ is
$f$-ample). Then $\Gamma = P$ and $f$ is isomorphic over $X$ to the
$(1,a,r-a)$ weighted blowup of $P \in X$.
\end{thm}

\begin{lemma}[Kawamata, \cite{Ka96}] \label{lem:Ka}
Let $P \in X \iso \frac{1}{r}(1,a,r-a)$ be as in Theorem \ref{thm:Ka}
and $f \colon (E \subset Y) \to (P \in X)$ the $(1,a,r-a)$ weighted
blowup; let $g \colon \widetilde{X} \to X$ be a resolution of
singularities with exceptional divisors $\{E_i \}$. Fix an effective
Weil divisor $H$ on $X$ and define $a_i =
a_{E_i}\left(K_X\right)$ and $m_i = m_{E_i}(H)$ in the usual way via
\begin{eqnarray}
 K_{\widetilde{X}} & = & g^*K_X + \textstyle\sum a_iE_i, \nonumber \\
 g^{-1}_*H         & = & g^*H - \textstyle\sum m_iE_i; \nonumber
\end{eqnarray}
define $a_E$ and $m_E$ similarly using $f$. Then
$m_i / a_i \le m_E / a_E$ for all $i$.
\end{lemma}

In \cite{Ka96} Lemma \ref{lem:Ka} is used to prove Theorem
\ref{thm:Ka}, but it is an interesting result in its own right; in
particular, it has two corollaries that are of great importance for
our problem.

\begin{cor} \label{cor:Kalemma:1}
Let $P \in X \iso \frac{1}{r}(1,a,r-a)$ and $f \colon (E \subset Y)
\to (P \in X)$ the Kawamata blowup as in Lemma \ref{lem:Ka}. Suppose
$\HH$ is a mobile linear system on $X$ and $n \in \ZZ_{\ge1}$ is such
that $P \in \CSXnH$. Then the valuation $v_E$ of $E$ is strictly
canonical or worse for $\XnH$.
\end{cor}

\begin{proof}[\textsc{Proof}]
Let $g \colon \widetilde{X} \to X$ be any
resolution of singularities with exceptional divisors $\{E_i\}$ and $H
\in \HH$ a general element. The assumption $P \in \CSXnH$ means that
$n \le m_i/a_i$ for some $i$, so by Lemma \ref{lem:Ka} $n \le m_E /
a_E$ as well.
\end{proof}

This Corollary \ref{cor:Kalemma:1} tells us that we
can exclude a singular point from any $\CSXnH$ simply by excluding the
valuation $v_E$; it has other uses as well.

\begin{cor} \label{cor:Kalemma:2}
Let $P \in X \iso \frac{1}{r}(1,a,r-a)$ and $f \colon Y
\to X$ be the Kawamata blowup of $P$ as in Lemma
\ref{lem:Ka}. Suppose $C \subset X$ is a curve containing $P$, $\HH$
is a mobile linear system on $X$ and $n \in \ZZ_{\ge1}$ is such that
$C \in \CSXnH$. Then $P \in \CSXnH$ also.
\end{cor}

\begin{proof}[\textsc{Proof}]
Let $g \colon \widetilde{X} \to X$ be a
resolution of singularities with exceptional divisors $\{E_i\}$ at
least one of which has centre $C$ on $X$ and is strictly canonical or
worse for $\XnH$. The rest of the proof is the same as that of
Corollary~\ref{cor:Kalemma:1}.
\end{proof}

\begin{absexclsingpts}
Suppose $P$ is a singular point of a hypersurface $X$ in one of
the~95 families and $P \in X$ is locally isomorphic to
$\frac{1}{r}(1,a,r-a)$. Let $f \colon (E \subset Y) \to (P \in X)$ be
the Kawamata blowup and suppose $B^3 < 0$ (where as always $B =
{-K_Y}$). We denote by $S$ the surface $f^{-1}_*\{x_0=0\}$,
which is an element of $\left|B\right|$ and is irreducible, assuming
$X$ is general.
\end{absexclsingpts}

\begin{lemma}[{see~\cite[5.4.3]{CPR}}] \label{Texists}
If $B^3 < 0$ then there exist integers $b,c$ with $b > 0$ and $b/r \ge
c \ge 0$ and a surface $T \in \left|bB + cE \right|$ such that
\begin{itemize}
\item[(a)] the scheme theoretic complete intersection $\Gamma = S \cap
T$ is a reduced, irreducible curve and
\item[(b)] $T\Gamma \le 0$.
\end{itemize}
\end{lemma}

\begin{thm} \label{thm:Tmethodabs}
Suppose $B^3 < 0$ and the integer $c$ provided by Lemma \ref{Texists}
is strictly positive. Then $P$ is excluded absolutely, that is, there
does not exist a mobile linear system $\HH$ of degree $n$ on $X$ such
that $K_X + \frac{1}{n}\HH$ is canonical and $P \in \CSXnH$.
\end{thm}

For the proof we need only the following two lemmas.

\begin{lemma} \label{lem:test_class}
Let $X$ be a Fano 3-fold hypersurface in one of the 95 families and
$\HH$ a mobile linear system of degree $n$ on $X$ such that $K_X +
\frac{1}{n}\HH$ is strictly canonical; suppose $\Gamma\subset X$ is an
irreducible curve or a closed point satisfying $\Gamma\in\CSXnH$, and
furthermore that there is a Mori extremal divisorial contraction
\[ f\colon (E\subset Y) \to (\Gamma\subset X), \quad \Centre_X E =
\Gamma, \]
such that $E\in\VXnH$. Then $B^2 \in \NEbar Y$.
\end{lemma}

\begin{proof}[\textsc{Proof}]
We know that
\[ K_Y + \textstyle\frac{1}{n}\HH_Y \qeq f^*\left(K_X +
\textstyle\frac{1}{n}\HH\right) \qeq 0 . \]
It follows that $B \qeq \frac{1}{n}\HH_Y$, and therefore $B^2 \in
\NEbar Y$, because $\HH_Y$ is mobile.
\end{proof}

\begin{lemma}[{see~\cite[5.4.6]{CPR}}] \label{NEbarY}
If $B^3 < 0$, let $T$ and $\Gamma = S \cap T$ be as in the conclusion
of Lemma \ref{Texists}. Write $R$ for the extremal ray of $\NEbar Y$
contracted by $f \colon Y \to X$ and let $Q \subset \NEbar Y$ be the
other ray. Then $Q = \RR_{\ge0}[\Gamma]$.
\end{lemma}

\begin{proof}[\textsc{Proof of Theorem \ref{thm:Tmethodabs}}]
Suppose $\HH$ is a mobile linear system of degree $n$ on $X$ such
that $K_X + \frac{1}{n}\HH$ is canonical and $P
\in \CSXnH$. Corollary~\ref{cor:Kalemma:1} of Kawamata's Lemma tells us
that the Kawamata blowup
$f \colon Y \to X$ of $P$
extracts a valuation $v_E$ (where $E = \Exc f$) which is strictly
canonical for $\XnH$. The test class Lemma~\ref{lem:test_class} now
implies that $B^2 \in \NEbar Y$.

But the ray $Q \subset \NEbar Y$ is generated by $(bB+cE)B$ for some
$b,c > 0$, and certainly $EB \in \NEbar Y$, so if $B^2 \in \NEbar Y$
we have both $B^2,EB \in Q$ (by definition of `extremal'). It follows
that $EB$ is numerically equivalent to $\alpha B^2$ for some positive
$\alpha \in \QQ$; but
\[ EB \cdot B = E \Big(A-\frac{1}{r}E \Big)^2 = \frac{1}{r^2}E^3 =
\frac{1}{a(r-a)} > 0,  \]
while
$B^2 \cdot B = B^3 < 0$
by assumption --- contradiction.
\end{proof}

\begin{cor} \label{cor:75singpts}
Let $X = X_{30} \subset \PP(1,4,5,6,15)_{x,y,z,t,u}$ be a general
member of family~75 of the~95 and $\HH$ a mobile linear system of
degree $n$ on $X$ with $\KXnH$ strictly canonical.  Then no singular
point of $X$ other than the two $\frac{1}{5}(1,4,1)$ points in the
$zu$-stratum can belong to $\CSXnH$.
\end{cor}

\begin{proof}[\textsc{Proof}]
Here is the complete list of singular points of $X$, together with the
sign of $B^3$ and $bB + cE \sim T$ (see Lemma~\ref{Texists}) for each
of them.
\begin{tabbing}
\hspace*{1cm} \=$P_y P_t \cap X = R_1, R_2 \sim
\frac{1}{2}(1,1,1)_{x,z,u}$ space \= sample text \= \kill
\>$P_y \sim \frac{1}{4}(1,1,3)_{x,z,u}$ \>$B^3 < 0$ \>$10B + E$ \\
\>$P_t P_u \cap X = P \sim \frac{1}{3}(1,1,2)_{x,y,z}$ \>$B^3 <
0$ \>$5B + E$ \\
\>$P_z P_u \cap X = Q_1, Q_2 \sim \frac{1}{5}(1,4,1)_{x,y,t}$
\>$B^3 < 0$ \>$4B$ \\
\>$P_y P_t \cap X = R_1, R_2 \sim \frac{1}{2}(1,1,1)_{x,z,u}$
\>$B^3 < 0$ \>$5B + 2E$
\end{tabbing}
Clearly all these apart from $Q_1$ and $Q_2$ satisfy the hypotheses of
Theorem~\ref{thm:Tmethodabs}, and consequently are excluded
absolutely.
\end{proof}

\begin{condexclsingpts}
Assume $P \in X$ is a singular point of a hypersurface in one of
the~95 families which is locally isomorphic to
$\frac{1}{r}(1,a,r-a)$ with $B^3 < 0$.  We keep the notations $T \sim
bB + cE$, $S = f^{-1}_*\{x_0=0\}$ and $\Gamma = S \cap T$ of
Lemma~\ref{Texists}; in the following paragraphs we consider the case
where the integer $c$ provided by Lemma \ref{Texists} is zero.  Out of
all such singular points, the vast majority live in a family with 
$b = a_1 < a_2$, where $a_1,a_2$ are the weights of $x_1,x_2$. For
such points we have the following result.
\end{condexclsingpts}

\begin{thm} \label{thm:Tmethodcond}
Let $P \in X$ be a singular point satisfying $B^3 < 0$, $T \sim bB$
and $b = a_1 < a_2$. Assume that the curve $C = \{x_0 = x_1 = 0\}\cap
X$ is irreducible. Then $R(Y,B) = k[x_0,x_1]$. It follows that if
$\HH$ is a mobile linear system of degree $n$ on $X$ such that $\KXnH$
is canonical and $P \in \CSXnH$ then in fact $\CSXnH =
\CS\left(X,\frac{1}{b}f_*\left|bB\right|\right)$, where $f \colon Y
\to X$ is the Kawamata blowup of $P$.
\end{thm}

\begin{proof}[\textsc{Proof for family 75}]
We do not prove this theorem here for every case.  As explained
in~\cite{Ry02}, the first statement, namely $R(Y, B) = k[x_0, x_1]$,
follows from two ray game calculations such as that in the proof of
our Lemma~\ref{lem:75tworay}; we have already shown this for
family~75.  The second part of the theorem follows easily from the
first: let $\HH$ be mobile of degree $n$ on $X$ with $\KXnH$ canonical
and $P \in \CSXnH$.  Then $\HH_Y \subset |nB| =
\left<x_0^n,x_0^{n-b}x_1,\ldots,x_0^rx_1^q\right>$, where $n = qb +
r$, $0 \le r < b$ --- so $\HH \subset f_*|nB| =
\left<x_0^n,\ldots,x_0^rx_1^q\right>$, while of course $f_*|bB| =
\left<x_0^b,x_1\right>$, and therefore $\CSXnH =
\CS\left(X,\frac{1}{b}f_*|bB|\right)$.
\end{proof}

We are now in a position to put together all the exclusion results
obtained so far to prove Theorem~\ref{thm:75aux} for family~75.

\begin{proof}[\textsc{Proof of Theorem~\ref{thm:75aux}}]
First, since $X$ is superrigid by~\cite[\S6]{CPR}, $\KXnH$ nonterminal
implies $\KXnH$ strictly canonical; therefore $\CSXnH$ is nonempty.
Corollary~\ref{cor:75curves} tells us that no curve belongs to
$\CSXnH$, Theorem~\ref{thm:smpts} says that the same is true for
smooth points and Corollary~\ref{cor:75singpts} shows the same for all
singular points other than $Q_1, Q_2 \sim
\frac{1}{5}(1,4,1)_{x,y,t}$.  Therefore at least one $Q_i \in \CSXnH$;
without loss of generality we may assume this holds for $Q_1$.

Now by Theorem~\ref{thm:Tmethodcond}
\[ \CSXnH = \CS\left( X, \textstyle\frac{1}{b} f_* |bB| \right) \]
where $b = 4$ because $T \sim 4B$ (in the notation above) and $f
\colon Y \to X$ is the Kawamata blowup of $Q_1$.  But $R(Y, B) =
k[x_0, x_1] = k[x, y]$, so $f_* |bB|$ is just $\left< x^4, y \right>$;
and both $x^4$ and $y$ have local vanishing order $4/5$ at $Q_2$, so
$Q_2 \in \CSXnH$ also, as required.
\end{proof}

\section{Birational classification of elliptic and K3
fibrations}  \label{sec:class}

We have now proved Theorem~\ref{thm:75aux} for family~75; all that
remains is to deduce the main Theorem~\ref{thm:75main} from it.  For
this we need the following theorem--definition and propositions.

\begin{thmdefn} \label{Kod}
Let $X$ be a variety, normal and
projective over $\CC$ as always, and $\HH$ a mobile linear system on
$X$. Fix $\alpha\in\QQ_{\ge0}$ and let $f \colon Y \to X$ be a birational
morphism such that $K_Y + \alpha\HH_Y$ is canonical. We define the
\emph{log Kodaira dimension $\kappa(X,\alpha\HH)$} to be
the $D$-dimension of $K_Y + \alpha\HH_Y$, that is,
\[ \kappa(X,\alpha\HH) = D(K_Y + \alpha\HH_Y) = \max
\{\dim(\im\fie_{\left|m(K_Y + \alpha\HH_Y)\right|}) \} , \]
taking the max over all $m\ge1$ such that $m(K_Y + \alpha\HH_Y)$ is
integral; if all the linear systems $\left|m(K_Y + \alpha\HH_Y)\right|$ are
empty, by definition \[\kappa(X,\alpha\HH) = D(K_Y + \alpha\HH_Y) =
{-\infty}.\] Then $\kappa(X,\alpha\HH)$ is independent of the choice of
$Y$ and is attained for a particular $Y$ using any sufficiently large
$m$ such that $m(K_Y + \alpha\HH_Y)$ is integral.
\end{thmdefn}

\begin{proof}[\textsc{Proof}]
This result is standard and is used in, e.g.,~\cite{Ch00}
and~\cite{Is01}.  The methods employed in \cite{Sh96} to show
uniqueness of a log canonical model can be used to prove it;
alternatively see \cite[Ch.\ 2]{FA}, particularly Theorem 2.22 and the
Negativity Lemma 2.19, which is an essential ingredient.
\end{proof}

\begin{emp}
Now let $X$ be a Mori Fano variety and $\HH$ a mobile linear system of
degree $n$ on $X$, that is, $n \in \QQ$ is such that $K_X +
\frac{1}{n}\HH \qeq 0$ --- of course $n\in\ZZ_{\ge 1}$ if $X$ is a
hypersurface in one of the 95 families.
\end{emp}

\begin{prop} \label{prop:varyep}
Assume that $K_X + \frac{1}{n}\HH$ is canonical.
\begin{itemize}
\item[(a)] Let $\ep \in \QQ$. Then
\[ \kappa\left(X, \left(\textstyle \frac{1}{n} + \ep \right)\HH \right) = 
\left\{\begin{array}{l@{\quad\mbox{if}\hspace{1.5ex}}l}
{- \infty} & \ep < 0 \\
0 & \ep = 0 \\
d \ge 1 & \ep > 0
\end{array} \right.\]
\item[(b)] If $1 \le \kappa(X, (\frac{1}{n} + \ep
)\HH) \le \dim X - 1$ for some $\ep \in \QQ_{>0}$ then $K_X +
\frac{1}{n}\HH$ is nonterminal, i.e., strictly canonical (so $K_X +
(\frac{1}{n} + \ep)\HH$ is noncanonical).
\end{itemize}
\end{prop}

\begin{proof}[\textsc{Proof}]
See the survey~\cite[III.2.3--2.4]{Is01}.
\end{proof}

Recall that the NFI-type inequality~\ref{nfi} was stated under two
alternative sets of assumptions: either
\begin{itemize}
\item[(a)]  $X$ is a Mori Fano and $\Phi \colon X \dashto Z/T$ a
birational map to the total space $Z$ of a $K$-trivial fibration $g
\colon Z \to T$; or
\item[(b)]  $X$ is a Mori Fano and $\Phi \colon X \dashto Z$ a
birational map to a Fano variety with canonical singularities.
\end{itemize}

\begin{prop} \label{dimT}
In situation (a) above, assume that $\KXnH$ is
canonical. Then for any $\ep \in \QQ_{>0}$, $\ka\XnepH = \dim T$.
\end{prop}

\begin{proof}[\textsc{Proof}]
Fix $\ep\in\QQ_{>0}$. $K_X + \frac{1}{n}\HH$ is canonical so
$\kappa(X,\frac{1}{n}\HH) = 0$, and therefore
$\kappa(Z,\frac{1}{n}\HH_Z) = 0$ by the birational invariance of log
Kodaira dimension (\ref{Kod}). But in fact $K_Z + \frac{1}{n}\HH_Z$ is
canonical (as is $K_Z + (\frac{1}{n} + \ep)\HH_Z$, because $K_Z$ is
canonical and $\HH_Z$ is free), so $\kappa(Z,\frac{1}{n}\HH_Z)$ (and
$\kappa(Z,(\frac{1}{n} + \ep)\HH_Z)$) can be computed on $Z$ as the
ordinary $D$-dimension. Consequently for $m \gg 0$ and such that $m(K_Z
+ \frac{1}{n}\HH_Z)$ is integral, it is in fact effective and
fixed. Fix such an $m$ with the additional property that $m\ep \in
\NN$, so that $m(K_Z + (\frac{1}{n}+\ep)\HH_Z)$ is integral as well;
let $F \in \left|m(K_Z + \frac{1}{n}\HH_Z)\right|$ be the unique
element. Now for any curve $C$ contracted by $g$, $FC = 0$, because by
assumption $K_Z C = 0$ and $\HH_Z = g^*\left|H\right|$. But $F$ is
effective, so it must be a pullback $g^* F_T$ of some effective $F_T$
on $T$. Furthermore, for any $m' \in \NN$,
\[ H^0(T,m'F_T) = H^0(Z,g^*(m'F_T)) = H^0(Z,m'F),  \]
so $F_T$ is fixed, $D(F_T) = 0$ and it is easy to see that $D(F_T +
m\ep H) = \dim T$, because $H$ is ample and $F_T$ is effective. Now
$g^*(F_T + m\ep H) = F + m\ep\HH_Z$, so
\begin{eqnarray*} \textstyle \kappa\left(X,\left(\frac{1}{n} + \ep
\right)\HH \right) & = & \textstyle D\left(m\left(K_Z +
\left(\frac{1}{n} + \ep \right)\HH_Z\right) \right) \\ & = & D(F+
m\ep\HH_Z) = \dim T
\end{eqnarray*}
as required.
\end{proof}

\medskip
\begin{proof}[\textsc{Proof of NFI-type inequality \ref{nfi} in
situation} \emph{(a)}]
We note that under the assumptions (a),~\ref{nfi} is an immediate
consequence of Propositions~\ref{prop:varyep} and~\ref{dimT}.  Under
assumptions~(b), the proof of~\cite[4.2]{Co95} can be easily adapted
to give an argument.  Like Theorem--Definition~\ref{Kod}, this is a
standard result used in~\cite{Ch00} and~\cite{Is01}, so we omit the
details.
\end{proof}

\bigskip
All that remains is to prove the main theorem for family~75.
Arguments similar to the following also prove
Theorems~\ref{thm:34main},~\ref{thm:88main} and~\ref{thm:90main}.

\begin{proof}[\textsc{Proof of Theorem \ref{thm:75main}}]
It is simplest to prove part (b) first; we then indicate how the
argument can be easily adapted to demonstrate (a) and (c) also.

{(b) \sloppy Suppose that $\Phi \colon X \dashto Z/T$ is a birational
map from $X$ to an elliptic fibration \linebreak $g \colon Z \to
T$. By the NFI-type inequality~\ref{nfi}, $\KXnH$ is nonterminal,
where as usual the system \linebreak[4] $\HH = \Phi^{-1}_*\HH_Z =
\Phi^{-1}_*g^*|H|$ is the transform of a very ample complete system of
Cartier divisors on $T$.  By Theorem~\ref{thm:75aux}, $\KXnH$ is
strictly canonical and $\CSXnH = \{Q_1,Q_2\}$. }

Let $Q$ be either $Q_1$ or $Q_2$. As in Lemma~\ref{lem:75tworay} we blow
up $Q$ and play out the two ray game; in the notation of the lemma,
\[ R(Y,B) = R(Y',B') = k[x,y], \]
$f'$ is the anticanonical morphism of $Y'$ and the composite $\pi
\colon X \dashto \PP^1$ is $(x^4,y)$. Now because $Q \in \CSXnH$, we
have that
\[ \HH_Y \subset |{-nK_Y}| = |nB| = k[x,y]_n =
k\left[x^n,x^{n-4}y,\ldots,x^{n \bmod 4}y^{\lfloor
n/4\rfloor}\right]; \]
the same is true of $Y'$ (since $Y$ and $Y'$ are isomorphic in
codimension one) and therefore \linebreak $\HH_{Y'} =
(f')^*\HH_{\PP^1}$ is the pullback of a mobile system on $\PP^1$. But
any mobile system on $\PP^1$ is free, so we can deduce (as in the proof of
Proposition~\ref{dimT}) that for any $\ep \in \QQ$ with $0 < \ep \ll 1$,
\begin{equation} \label{kale1}
\ka\XnepH = \ka\left(Y',\left(\textstyle\frac{1}{n} +
\ep\right)\HH_{Y'}\right) = D(F + m\ep\HH_{\PP^1}) \le 1
\end{equation}
where $F$ is a fixed effective divisor on $\PP^1$ (so in fact $F = 0$)
and $m \in \ZZ_{>0}$.

But by Proposition~\ref{dimT} applied to $\Phi \colon X \dashto Z/T$
we have $\ka\XnepH = \dim T = 2$, which contradicts
(\ref{kale1}). This proves (b).

(a) We can follow the proof of (b) but in the end, rather than a
contradiction, we deduce that the system $\HH = \Phi^{-1}_*g^*|H|$ is
actually a pullback $\pi^{-1}_*\HH_{\PP^1}$ via the map $\pi = (x^4,y)
\colon X \dashto \PP^1$. This induces an isomorphism $\PP^1 \to T$
such that the specified diagram commutes.

(c) If we assume $\Phi$ is not an isomorphism, we can follow the
argument for (b) to deduce that $\ka\XnepH = 1$, which is obviously a
contradiction since $\HH_Z = |H|$ is very ample. (Note that the
NFI-type inequality~\ref{nfi} requires us to assume $\Phi$ is not an
isomorphism in the Fano case; for the elliptic and K3 cases this is of
course not necessary.)

This completes the proof.
\end{proof}

\end{document}